\documentclass[11pt]{article}
\usepackage{latexsym,amsmath}
\usepackage{graphicx}
  \ifx\LabelFigloaded\MYundefined\relax
  \else
   \fi

  \def\LabelFigloaded{\relax}

  \chardef\LabelFigCatAt\the\catcode`\@
  \catcode`\@=11

 \let\LabelFigwlog@ld\wlog
 \def\wlog#1{\relax}

 \ifx\\\MYundefined@
    \let\\\relax
 \fi

  \def\ms@g{\immediate\write16}

 \def\N@wif{\csname newif\endcsname }
 \def\Temp@ {\N@wif\ifIN@}
 \ifx\INN@\MYundefined@
    \else \let\Temp@\relax
 \fi
 \Temp@

  \def\IN@{\expandafter\INN@\expandafter}
  \long\def\INN@0#1@#2@{\long\def\NI@##1#1##2##3\ENDNI@
    {\ifx\m@rker##2\IN@false\else\IN@true\fi}%
     \expandafter\NI@#2@@#1\m@rker\ENDNI@}
  \def\m@rker{\m@@rker}
 
  \newtoks\Initialtoks@  \newtoks\Terminaltoks@
  \def\SPLIT@{\expandafter\SPLITT@\expandafter}
  \def\SPLITT@0#1@#2@{\def\TTILPS@##1#1##2@{%
     \Initialtoks@{##1}\Terminaltoks@{##2}}\expandafter\TTILPS@#2@}

 \def\Shifted@@#1#2#3{\setbox0=\hbox{#3}%
   \raise -\dp0\vbox {\kern-#2%
       \hbox {\kern#1\unhbox0\kern-#1}%
           \kern#2}}

 \newcount\gridcount
 \newbox\auxGridbox@ \newbox\hGridbox@ \newbox\vGridbox@
 \newbox\Labelbox@ \newbox\auxLabelbox@
 \newbox\Coordinatebox@
 \newtoks\Labeltoks@
 \newdimen\Wdd@ \newdimen\Htt@
 \newdimen\Wddd@ \newdimen\Httt@
 
 \def\Wr@{\immediate\write16}

 \newdimen\GL@wd
 \def\GridLineWidth#1{\GL@wd=#1}

 \def\gobble#1{}
 \def\EdgeErr@{\Wr@{}%
      \Wr@{\string\Edges\space argument
      1, 10, 100 or 1000 please\string!}%
      }

 \newcount\Edgect@

 \def\Sweepup#1\endSweepup{}

 \def\SetEdges@{%
    \edef\Zr@@s{\expandafter\gobble\number\Edgect@\empty}%
        \count255=0\Zr@@s\relax
        \ifnum\count255=\z@\else\EdgeErr@\show\tailtest\fi
        \count255=1\Zr@@s\relax
        \ifnum\count255=\Edgect@\relax\else\EdgeErr@\show\leadtest\fi
    \EdgGl@b\edef\Zr@s{\expandafter\gobble\Zr@@s\empty}
    \ifnum\Edgect@>\@ne\relax\EdgGl@b\let\L@Dc\empty
        \else\EdgGl@b\edef\L@Dc{\string.}\fi
    \ifnum\Edgect@>\@ne\relax
        \EdgGl@b\edef\Edgescale@##1{\divide##1 by \Edgect@}%
        \else\EdgGl@b\edef\Edgescale@##1{}\fi
    }

 \def\Edges#1{\Edgect@=#1\relax
     \let\EdgGl@b\global \SetEdges@}

 \Edges{1}

 \def\hhrule{\hrule height \GL@wd\vskip-.\GL@wd}

 \def\hRule@{%
   \advance\gridcount -2%
   \vfil\hhrule\vfil
   \llap{\smash{\raise -2.5pt
     \hbox{\L@Dc\number\gridcount\Zr@s\kern2pt}}}%
   \hhrule
   }

\def\vvrule{\vrule width \GL@wd \kern-\GL@wd}

 \def\vRule@{\advance\gridcount 2%
   \hfil\vvrule\hfil
   \setbox\auxGridbox@=\vbox to 0pt
      {\vskip \Htt@\vskip 2pt
        \hbox to 0pt{\hss\L@Dc\number\gridcount\Zr@s\hss}\vss}%
      \wd\auxGridbox@=0pt \box\auxGridbox@
   \vvrule
   }

 \def\PlaceGrid@@{\gridcount=10 
  \setbox\hGridbox@=\hbox{%
        \hbox{%
             \hskip-.4pt\vrule
             \vbox to \Htt@{%
               \offinterlineskip\parindent=\z@\relax
               \hbox to \Wdd@{\hfil}
               \hRule@\hRule@\hRule@\hRule@
               \vfil\hhrule\vfil}%
             \vrule\hskip-.4pt}
    }%
  \gridcount=0%
  \setbox\vGridbox@=\hbox{%
      \vbox{\offinterlineskip\parindent=0pt\hsize=0pt
         \vskip-.4pt\hrule%
         \hbox to \Wdd@{%
                 \vtop to \Htt@{\vfil}%
                 \vRule@\vRule@\vRule@\vRule@
                 \hfil\vvrule\hfil}%
         \hrule\vskip-.4pt}}%
  \wd\hGridbox@=0pt\ht\hGridbox@=0pt
  \wd\vGridbox@=0pt\ht\vGridbox@=0pt
  \hbox{\box\hGridbox@\box\vGridbox@}%
  }

 \def\LabelsGlobal{\def\LabGl@b{\global}}
 \def\LabelsLocal{\def\LabGl@b{}}
 \LabelsGlobal 

 \def\SetLabels#1\endSetLabels{%
   \LabGl@b\Labeltoks@={#1()\\}%
   }

 \LabGl@b\Labeltoks@={()\\}

 \def\ShowGrid{\LabGl@b\let\PlaceGrid@\PlaceGrid@@}
 \def\HideGrid{\LabGl@b\let\PlaceGrid@\relax}
 \def\Grids{\ShowGrid\LabGl@b\let\GridSwitch@\ShowGrid}
 \def\noGrids{\HideGrid\LabGl@b\let\GridSwitch@\HideGrid}

 \noGrids

 \def\bAdjust@@{%
     \setbox\auxLabelbox@=\hbox{\raise \dp\auxLabelbox@
            \box\auxLabelbox@}}
 \def\bAdjust@{\let\vAdjust@\bAdjust@@}

 \def\eAdjust@@{\dimen0=-.5\ht\auxLabelbox@
     \advance\dimen0 by .5\dp\auxLabelbox@
     \setbox\auxLabelbox@=
            \hbox{\raise\dimen0\box\auxLabelbox@}}
 \def\eAdjust@{\let\vAdjust@\eAdjust@@}

 \def\tAdjust@@{%
     \setbox\auxLabelbox@=\hbox{\raise-\ht\auxLabelbox@
            \box\auxLabelbox@}}
 \def\tAdjust@{\let\vAdjust@\tAdjust@@}

 \let\vAdjust@\relax

 \def\lAdjust@{\let\hAdjust@\rlap}
 \def\rAdjust@{\let\hAdjust@\llap}

 \let\hAdjust@\relax\let\vAdjust@\relax

 \def\FetchLabel@#1(#2)#3\\{%
     \IN@0#2@@\ifIN@
        \setbox0=\hbox{\ignorespaces#1#3\unskip}%
        \ifdim\wd0>0pt
           \ms@g{}%
           \ms@g{ !!! Bad label(s)? !!!}%
           \message{ #1(#2)#3}%
        \fi
        \def\LabelMole@##1\endFetchLabel@{%
            \IN@0()\\@##1@%
            \ifIN@\def\Temp@{\FetchLabel@##1\endFetchLabel@}%
            \else\def\Temp@{}%
            \fi
            \Temp@
           }%
     \else
       \ignorespaces#1\unskip
       \setbox\auxLabelbox@=%
         \hbox to 0pt{\hss\ignorespaces\hAdjust@
          {\ignorespaces#3\unskip}\hss}%
       \vAdjust@
       \let\hAdjust@\relax\let\vAdjust@\relax
       \AugmentLabelBox@@{#2}%
       \ht\Labelbox@=0pt\dp\Labelbox@=0pt
       \let\LabelMole@\FetchLabel@%
     \fi\LabelMole@}

 \newtoks\XYSep@ 
 \def\SetXYSeparator#1{%
     \IN@0#1@@\ifIN@\XYSep@{*}%
     \else
     \XYSep@{#1}%
     \fi
     }

 \SetXYSeparator*

 \def\AugmentLabelBox@@#1{%
     \IN@0\the\XYSep@ @#1@\ifIN@
       \SPLIT@0\the\XYSep@ @#1@%
       \setbox\Labelbox@=\hbox to 0pt{%
         \unhbox\Labelbox@
         \Shifted@@{\the\Initialtoks@\Wddd@}%
         {\the\Terminaltoks@\Httt@}%
         {\box\auxLabelbox@}}%
     \else
         \ms@g{}%
         \ms@g{ !!! Bad insertion point. !!!}%
         \message{ (#1\ this point was rejected.)}%
     \fi
    }

 \def\FetchOption@#1[#2]#3\endFetchOption@{%
    \def\temp{#1}
    \ifx\temp\empty
       \Edgect@=#2\relax
       \let\EdgGl@b\relax
       \SetEdges@
       \Cleaner@#3%
    \fi}

 \def\Cleaner@#1[@]{\Labeltoks@{#1}}
     
 \def\PlaceLabels@@{\mathsurround=0pt
     \def\Cr@{\\}%
     \let\L\lAdjust@\let\R\rAdjust@
     \let\B\bAdjust@\let\E\eAdjust@\let\T\tAdjust@
     \expandafter\FetchOption@\the\Labeltoks@[@]\endFetchOption@
     \Wddd@=\Wdd@ \Edgescale@\Wddd@ 
     \Httt@=\Htt@ \Edgescale@\Httt@
     \expandafter\FetchLabel@\the\Labeltoks@\endFetchLabel@
     \box\Labelbox@
     }%

 \let \PlaceLabels@\PlaceLabels@@

 \def\AffixLabels#1{\setbox\Coordinatebox@=\hbox{#1}%
      \Wdd@=\wd\Coordinatebox@ \Htt@=\ht\Coordinatebox@
      \advance\Htt@ \dp\Coordinatebox@
      \hbox{\copy\Coordinatebox@\kern-\Wdd@ 
           \Shifted@@{0pt}{-\dp\Coordinatebox@}%
           {\PlaceLabels@\PlaceGrid@}%
           \kern\Wdd@}%
      \GridSwitch@ 
      \LabGl@b\Labeltoks@{()\\}%
      }
 
   \let\wlog\LabelFigwlog@ld   
   \catcode`\@=\LabelFigCatAt  

                                By

              Raymond S\'eroul <A18645@FRCCSC21.BITNET>
                                and 
              Laurent Siebenmann <lcs@topo.math.u-psud.fr>
    
              VERSIONS: July 1991, Oct 1991, Jan 1992, July 1992

INTRODUCTION

      This labelling package is intended for TeX users who
rely on non-TeX sources for for their graphics inserts.  It
provides means for adding TeX labels to such inserts with a
minimum of fuss. 

       For most labels, TeX users have in the past found it
reasonably convenient to rely on non-TeX sources. Typical
occasions when an inescapable need for TeX labels seemed to
arise are

 (a) when the graphics program lacks certain exotic or complex
mathematical symbols

 (b) when the very highest typographical quality is wanted for the
labels

 (c) when labels included with the graphics fail to print, 
 and you cannot figure out why (cf. boxedeps.doc).  The labels
 provided by labelfig.tex are 100

       Since this package first appeared, many users, who in the
past scarcely dreamed of using TeX labels, have come to use
nothing but.  So it is now appropriate to add

Intoxication Warning:  TeX labels may be addictive and expensive. 

     If you have a fast preview you may disagree, and even find
that this package provides an agreeable paste-up environment; see
extra applications at end.

     Note to publishers: It is possible and convenient to ultimately
export the TeX labels produced by labelfig.tex to become an integral
part of the EPS file. This is often desired by a publisher who typically
uses an "upmarket" graphics or page layout program, with which the
staff is skilled in perfecting figures.  See Appendix I for
a recipe.

     The authors are grateful to Patrick Ion of Math Reviews for
helpful comments and encouragement.

BASIC INSTRUCTIONS

    After reading in the macro file using

preview or proof your figure with a coordinate grid printed on
top, by typing the following:

    \ShowGrid  
    \AffixLabels{<the graphics insertion>}

Here <the graphics insertion> is what you would type to insert
the graphics object alone without the grid.  This must provide
for the space around it. For example <the graphics insertion>
might well be \BoxedEPSF{MyFigure scaled 700} using the
boxedeps.tex macro package (from same source); this provides a
TeX box containing the encapsulated PostScript insert specified by
the file MyFigure. \AffixLabels{...} provides the grid (supposing
\ShowGrid is present) and later, once you have specified labels
using the grid, it will "tack on" the labels.

     The grid is a sort of (usually elongated) checkerboard of
ten rows and ten columns and its (internal) partitions are by
default numbered  .1, ... ,.9  both horizontally (X-coordinate
running left to right) and vertically (Y-coordinate running bottom
to top).  Thus the points enclosed by the grid correspond to the
points of the unit square in the cartesian "X-Y" plane, the lower
left corner corresponding to the origin (0,0).  By extrapolation,
the full page corresponds to a larger rectangle in the plane.

     These coordinates serve to position labels as follows.
Before the \AffixLabels{...} command type label specifications:

  \SetLabels
   (<X-coordinate>*<Y-coordinate>) <first label> \\
   .
   .
   .
   (<X-coordinate>*<Y-coordinate>)  <last label> \\
  \endSetLabels

Each row specifies one label and is terminated by \\.  In each
row, the position indicator comes first; it is written as a
standard cartesian point except that the X- and Y- coordinates
are separated by * rather than a comma because TeX allows a
comma as decimal point. There are no dimension units to specify
as the unit is the grid itself.

     By default, this cartesian point specifies where the middle
of the baseline of the label will be located.  However if you precede
the point by \L [or \R] the left [or right] edge of the baseline will
be located there. Similarly you may also precede the point by \T, \E,
or \B to vertically align the top equator or bottom of the label box
at the specified point.  This gives nine standard positions of
the label with respect to the insertion point --- corresponding to
the eight principle points of the compas and the center

                     \L\T     \T      \R\T

                     \L\E     \E      \R\E

                     \L\B     \B      \R\B

But this neglects the default "baseline" level of TeX,
giving potentially three more positions

                     \L    <no tag>   \R

For text, the baseline level is often the preferred. Its relation to
the others is variable. It will often coincide with the bottom level,
as happens for "X".  But it is often distinct, as for "g", in which
case you have in all 12 distinct positions rather than 9.

     It is convenient to think of this specification of label
position as attaching the label by a thumb-tack to the coordinate
grid. There are up to twelve positions of the thumb-tack on the
label, while the position of the thumb-tack on the coordinate grid is
arbitrary.  Normally, one choses the position of the thumb-tack on
the label to be the one that is the closest to the item being
labeled.  There are good reasons for this "rule of thumb":

   (a)  It facilitates correct positioning at first try.

   (b)  If the scale of the figure must be altered after labels
have been affixed, the labels have a good chance of remaining well
positioned.

   (c)  The visible grid need not extend beyond the "bounding box"
for the figure, because the best preferred position is always
(at least almost) within the bounding box .

The second reason is particularly important. Indeed it often
happens that scale has to be altered after labelling begins, in
order to either provide space for the labels, or to adjust
proportions between the labels and the figure.  (The size of labels
is unaffected by scaling.)

     Here is an artificial but self-contained test which uses
TeX rules to make a graphics object.

TEST

    Do not skip this!

 \def\FrameIt#1{\hbox{\vrule$\vcenter {\hrule\kern3pt%
             \hbox {\kern3pt #1\kern3pt}%
               \kern3pt\hrule}$\relax\vrule}}

 \def\Caption#1#2{\FrameIt{%
       \vtop {\hsize=#1\relax \parindent=0pt
         \leftskip=0pt \rightskip=0pt plus15pt
         \parfillskip=0pt
         \lineskip=1pt\baselineskip=0pt
         #2}}}

 \def\FirstQuadrant{\hbox to 100pt{\vrule\vbox to 100pt{%
        \hbox to 100pt{\hfil}\vfil\hrule}\hss}}

  \SetLabels
    \R(.5*.2) $\zeta\,\cdot$\\
    (.9*-.10) $\xi$\\
    \R(-.03*.9) $\eta$\\
    \T(.5*.9) \Caption{70pt}{%
          \it The norm of
          $g(\xi+i\eta)$ is indicated on
          contours of this invisible surface.}\\
  \endSetLabels

  \AffixLabels{\FirstQuadrant}

  \end

  Note that the coordinates to use for labels are indicated on the
edges of the grid (when visible) corresponding to the conventional
x- and y- axes of the Cartesian plane. By default the grid is
1-by-1. However, by the command \Edges{100}, you can change this
to 100-by-100 and many users find this alternative most
convenient. Place the command \Edges{...} in your style file (or
header) since its effect is is global. Other possible edge values
are 10 and 1000.

  If you use the command \Edges{...} at all, do so with care.  For
if you accidentally delete an \Edges{...} command your labels will
abruptly be badly misplaced and may logically but mysteriously
generate "dimension too big" errors under TeX and "off page" errors
under your driver.  

  You can dictate the edgescale for an individual figure by giving
the scale in brackets immediately after \SetLabels.  Thus, to
import into an article using say \Edge{100} a figure labelled using
another edgescale, say the original 1-by-1 default, you can use
\SetLabels[1]...\endSetLabels.

GETTING IT DOWN PAT

     Complicated labeling deserves the same respect as
complicated mathematics.  Do not expect it to come out perfect the
first time!  What is needed in either case is a mechanism to
repeatedly typeset troublesome pieces.

     One mechanism is always available.  One does complicated
labelling in a separate "test" file involving just the figure being
labelled;  a texpert will know how to \dump TeX's current state as
a temporary format that restarts rapidly at each retry.  Usually,
one then pastes the completed labelled figure back into the main
TeX file, but, of course, one can also \input it as an auxiliary
file.

     If you do not have a TeXpert at handy, here is a first
approximation to an efficient setup. By deletions reduce a copy
of your article to just a few lines before and after the figure.
Now label the figure, and finally, copy and paste the labelled
figure to the original article. Then copy the next figure to label
into this testbed and repeat. The TeXpert can improve the  speed
at which TeX starts up, by compiling a format specifically for
your article; just one caution: best NOT include in the format
ephemeral details of setup like \Set<mydriver>ArtSpecials (from
boxedeps.tex because this reads  figure dimensions which you may
change during your work session.

     An improved mechanism to repeatedly typeset troublesome
pieces is now available on the Macintosh; it is called LinoTeX;
see the same ftp sources.  It could be set up on many types
of computer.

     Before using labelfig.tex to attach labels to a graphics
object inserted using boxedeps.tex or BoxedArt.tex, make it a
firm rule to carefully adjust the bounding box using the trimming
commands of these packages, and also at least tentatively scale
and position the object. Beware of changing the grid inadvertently
after the labels have been positioned.  For example, correcting
the bounding box of a PostScript graphics object can foul up the
labels by changing the coordinate grid to which the labels are
attached. This is particularly true for the trimming  commands of
boxedeps.tex and BoxedArt.tex. However, as noted already, change
of scale is much less disruptive, and modest adjustments should be
well tolerated.

     Sometimes the labels protrude so far from the bounding box
of a figure that the figure has to be repositioned.  Best do this
by ad hoc spacing, say using \hglue and \vglue; altering the
bounding box would create a vicious circle.

     Remember that you are responsible for preventing labels
from overlapping. You are responsible for all label typography
including size and style. A label is really just about anything
that can be put in a TeX box. Note that spaces at the beginning
and end of labels will normally be suppressed; if you really want
them you must protect them with TeX braces.

     This package temporarily sets the \mathsurround parameter
of TeX to zero  while the labels are being affixed. This is done
because nonzero \mathsurround space would influence the position
of left and right aligned labels; then, when a texpert or printer
modifies mathsurround, diagram labeling might be disastrously
altered. There is a small price to pay involving labels that are
formatted as caption boxes including mathematics: you  may want or
need to specify an explicit mathsurround space within the caption
box; it will not influence anything outside.

     Those hostile to the use of * as separator between
the X and Y coordinates of label insertion points, are free to
impose another using \SetXYSeparator{<the new separator>}.  
Americans may prefer "," to "*" since they never use a 
comma as a decimal point; on the other hand, * may be more visible.

APPENDIX (I)  MERGING labelfig.tex LABELS INTO AN EPSF GRAPHICS OBJECT.

     As promised in the introduction, here is a recipe useful for
publishers. It works at least on Macintosh and at least for vectorized
graphics and Adobe type1 fonts.  (There is surely a similar recipe for
PCs under MSWindows.)

 (a)  Use boxedeps.tex utility to integrate the figure given by the eps
file, "x.eps" say, with a visible frame around it.  See
\ShowDisplacementBoxes command in boxedeps.tex.  To get precise results
automatically it is important to use the \Trim... commands of
boxedeps.tex making the "DisplacementBox" neatly fit the figure.

 (b)  Use the TeX printer driver and LaserWriter (versions >= 8.1.1) to
export to an EPSF the DVI page containing the integrated, labelled
figure. You now have an EPS file  "xx.eps"  that contains too much, and at
the wrong scale, and at wrong position.

 (c)  Convert the EPSF to an Adode Illustrator format EPSF using
the shareware utility called epsConvert by Sam Weiss
1993-- (currently $25).

 (d)  In Illustrator (or a compatible program), group the labels and the
"DisplacementBox"; copy them to the clipboard and paste them into "x.ps".
This step requires that all the label fonts be "visible to the Macintosh.

 (e)  Translate and scale the pasted group consisting of the labels plus
the "DisplacementBox" so as to make the "DisplacementBox" the bounding
box of (labelless) figure represented by "x.eps".  At this point the
labels will be correctly placed on the figure "x.eps".

 (f)  Ungroup and delete the "DisplacementBox".  The result is the
desired single EPS file, "x+.eps" say, It contains the original figure
plus its labels.  

     Using grouping and ungrouping appropriately in "x+.eps", a
publisher's staff can very efficiently improve label positions etc.

APPENDIX II)  SOME EXOTIC APPLICATIONS

     The grid of labelfig.tex is analogous to a light-table in
classical page makeup with wax or latex glue.  In principle, you
can use it to compose any page from its indivisible parts.  This
even has some of the artisanal charm of classical paste-up
provided you have a fast screen preview to make the process
"interactive".

     In practice labelfig.tex is a tool for nonstandard jobs.
Here are a few going beyond the labelling already discussed.

(I)  GRAPHICS INTEGRATION.

     This is accomplished by treating the imported graphics
objects as labels.  The underlying graphics object is then
typically an empty  \vbox to <dimension>{\vfill} in a TeX
\midinsert...\endinsert construction.  A label line
might be of the form

   (.1*.1) \\

The exact form of the special command varies from driver to
driver.  However, in the case of encapsulated PostScript graphics
(EPSF norm), by relying on boxedeps.tex, one can have the
following standard syntax (independant of driver  (see
boxedeps.doc for details.
  
  (.1*.1) \BoxedEPSF{MyFigure scaled <scale in mils>}\\

This may be slow since it requires TeX to read the PostScript
file to read bounding box using many complex macros.  So you
may want to try

  (.1*.1) \EPSFSpecial{MyFigure}{<scale in mils>}\\

which is fast and driver independant, but it squashes the
bounding box, normally to its lower left corner.

     Similarly for graphics of the Macintosh PICT norm ---
using BoxedArt.tex (same sources) in place of boxedeps.tex.

     This approach to integration is to be recommended when
one is assembling a composite graphics object.

 (II)  COMMUTATIVE DIAGRAM ENHANCEMENT

     Commutative diagrams or arrays of mathematical objects
connected by arrows of various sorts are common in mathematics.
The mathematical objects require the use of TeX.  Recently TeX
acquired a good collection of arrows of all slopes --- that of
LamSTeX --- plus pwerful macros to build the diagrams.

     However, even the LamSTeX collection is often
inadequate; it lacks for example double shafted arrows, dotted
arrows and curved arrows. Fortunately it is possible to produce
such arrows on an individual basis using sophisticated graphics
programs such as Illustrator and AldusFreehand (both serving
the EPSF norm) or using Metafont (with its public domain norm).
Since the creation of each new arrow is a work of love, you
probably want to limit the number of arrows by using LamSTeX
for most arrows. The 40K commutative diagram module of LamSTeX
has been adapted to work with AmSTeX and a copy may be posted
with LabelFig and related files. Unfortunately no one has yet
offered a version that works with Plain TeX or LaTeX.

       Suffice it here to say that when the exotic arrow has
been somehow imported into TeX, labelfig.tex treats it as a
label that one affixes to the commutative diagram.  Two other
steps will be treated in separate notes, namely the matter of
extracting the dimension specifications for the arrow and the
construction of the arrow --- for these steps are far from
unique and often depend intimately on your computer environment. 
Notes for the Macintosh-Textures-Illustrator combination are
found in the file ExoticArrows.doc.

 (III) NESTING 

Ingenuity pays off in exploiting labelfig.tex. One can
mix graphics and typography quite freely.  labelfig.tex is good
for freeform or overlapping arrangements, while boxedeps.tex (or
BoxedArt.tex) is best for regimented non-overlapping
arrangements --- and the two can be combined.

     The default behavior of labelfig.tex is not ideal 
for nesting objects, because to prevent trouble for beginners
the register for labels is globally cleared when \AffixLabels
concludes.  But there are switches available

      \LabelsGlobal      \LabelsLocal

which change this.  To understand this, extend the above test 
by something like:

 \LabelsLocal

 \SetLabels
    (.5*.5) AAA\\
 \endSetLabels

 {
 \SetLabels
    (.5*.5) ZZZ\\
 \endSetLabels
   \AffixLabels{\FirstQuadrant}
 }

   \AffixLabels{\FirstQuadrant}

     There are however potential pitfalls.  Neither
labelfig.tex nor boxedeps.tex has been tested under extreme
conditions. Problems may occur if their procedures are
indiscriminately nested. For boxedeps.tex (not labelfig.tex)
there is a precise cause for worry, namely many of its
variables are "global", which means that TeX braces will not
provide the protection one might expect.

COMMAND SUMMARY FOR labelfig.tex

  Here [...] means optional (one or zero)
       [...]* means any number of such constructs

  \SetLabels
    [[<P>](<X><Sep><Y>) <label> \\]*
  \endSetLabels
  \ShowGrid  
  \AffixLabels{<the figure>}

   --- <P> is tack position, one of eleven or empty
              order irrelevant

                   \L\T      \T      \R\T

                   \L\E      \E      \R\E

                     \L               \R

                   \L\B      \B      \R\B

   --- (<X><Sep><Y>) insertion point;
  <Sep> is separator, = * by default;
  \SetXYSeparator{<Sep>} changes it.
   <X> and <Y> are real numbers

  --- <label> a label to attach 

  --- <the figure> the figure to label 

  \GlobalLabels (default)     
  \LocalLabels  setting for nested constructs.

 \Grids makes ALL grids appear; \HideGrid then makes just next disappear.
 \noGrids returns to default.  The commands are always global.

 \GridLineWidth{<dimension>} adjusts width of grid lines. Default is very
small, to give "hairline" effect. If your grid lines are missing try
setting \GridLineWidth{1pt}.

 \Edges#1 globally changes the edge size of all grids to the numerical 
value #1, which must be 1, 10, 100, or 1000.  The default is 1.

VERSION HISTORY.
 --- Jan 1993: \Edges#1 and [??] option after \SetLabels
 --- July 1992: \Grids, \noGrids, \HideGrid;
       Gridlines become hairlines; \GridLineWidth{<dimension>}.
 --- Oct 1991, Jan 1992: \SetXYSeparator{<Sep>},  \LabelsGlobal,
       \LabelsLocal.
 --- July 1991: first release

Address for bugs and other feedback:

        Raymond S\'eroul
        IREM and Lab. de Typographie Informatise
        Univ. Rene Descartes
        Strasbourg

    Tel 33-88-41-63-45
    Email:  A18645@FRCCSC21.BITNET

        Laurent Siebenmann
        Mathematique, Bat. 425,
        Univ de Paris-Sud,
        91405-Orsay,
        France

    Tel 33-1-6941-7949; 
    Email: lcs@topo.math.u-psud.fr

  \newcommand{\lab}[1]{\label{#1}}                

\def\fullpage {
\addtolength{\topmargin}{-1.5 cm}
\addtolength{\oddsidemargin}{-1.5 cm}
\addtolength{\textwidth}{+3 cm}
\addtolength{\textheight}{+3 cm} }

\fullpage

\usepackage{amssymb}
\usepackage{amsthm}
\usepackage{epsfig}
\begin{document}
\newtheorem{theorem}{Theorem}
\newtheorem{cor}{Corollary}
\newtheorem{lemma}{Lemma}
\newtheorem{prop}{Proposition}

\newcommand\eps{\varepsilon}
\newcommand{\E}{\mathbb E}
\newcommand{\Var}{{\rm Var}}
\newcommand{\Prob}{\mathbb{P}}
\newcommand{\N}{{\mathbb N}}
\def\lam{\lambda}
\def\la{\lambda}
\def\G{{\cal G}}
\def\pr{{\bf P}}
\def\ex{{\bf E}}
\renewcommand{\E}{\ex}
\renewcommand{\Prob}{\pr}

\newcommand{\be}{\begin{equation}}
\newcommand{\ee}{\end{equation}}
\newcommand{\bea}{\begin{eqnarray}}
\newcommand{\eea}{\end{eqnarray}}
\newcommand{\non}{\nonumber}
\newcommand{\bean}{\begin{eqnarray*}}
\newcommand{\eean}{\end{eqnarray*}}
\newcommand\eqn[1]{(\ref{#1})}
\newcommand{\bel}[1]{\be\lab{#1}}
\newcommand{\lbel}[1]{\be}

\def\blackslug{\hbox{\kern1pt\vrule height6pt width4pt  depth1pt\kern1pt}}
\def\qed{\penalty 500\hbox{\quad\blackslug}\ifmmode\else\par
  \vskip4.5pt plus3pt minus2pt\fi}

\def\no{\noindent}
\def\bs{\bigskip}
\def\ss{\smallskip}
\def\ms{\medskip}
\def\K{{\cal K}}
\def\G{{\cal G}}

\title{On the threshold for $k$-regular subgraphs \\ of\\  random graphs.}
\author{
Pawe{\l} Pra{\l}at\\
{\small Department of Mathematics and Statistics} \\
{\small Dalhousie University} \\
{\small Halifax NS, Canada} \\
\and Jacques Verstra\"{e}te\\
{\small Department of Mathematics and Statistics} \\
{\small McGill University}\\
{\small Montreal QC, Canada}
\and Nicholas Wormald
\thanks{Supported by the Canada Research Chairs Program and NSERC}\\
{\small Department of Combinatorics and Optimization} \\
{\small University of Waterloo}\\
{\small Waterloo ON, Canada}}
\date{}
\maketitle

\begin{abstract}
The $k$-core of a graph is the largest subgraph of minimum degree at
least $k$. We show that for $k$ sufficiently large, the $(k +
2)$-core of a random graph $\G(n,p)$ asymptotically almost surely
has a spanning $k$-regular subgraph. Thus the threshold for the
appearance of a $k$-regular subgraph of a random graph is at most
the threshold for the $(k+2)$-core. In particular, this pins down
the point of appearance of a $k$-regular subgraph in $\G(n,p)$ to a
window for $p$ of width roughly $2/n$ for large $n$ and moderately
large $k$.
\end{abstract}

\section{Introduction}

In this paper, we study the appearance of $k$-regular subgraphs of
random graphs. The {\sl $k$-core} of a graph $G$ is the unique
largest subgraph of $G$ of minimum degree at least $k$ (note that
the $k$-core may be empty). Evidently, the $k$-core of a graph can
be found be repeatedly deleting vertices of degree less than $k$
from the graph. In the case $k = 2$, this corresponds to the
appearance of cycles in $\G(n,p)$, which is well-researched, and
precise results concerning the distribution of cycles may be found
in Janson~\cite{J} and Flajolet, Knuth and Pittel~\cite{FKP}. By
analysing the vertex deletion algorithm for the Erd\H{o}s-R\'{e}nyi
model $\G(n,p)$ of random graphs, Pittel, Spencer and
Wormald~\cite{PSW} proved that for fixed $k \geq 3$, there exists a
constant $c_k$ such that $c_k/n$ is a sharp threshold for a $k$-core
in $\G(n,p)$. (When discussing thresholds of $k$-cores and
$k$-regular subgraphs,
 we mean thresholds for nonempty $k$-cores and nonempty $k$-regular subgraphs.) Here \bel{ckdef}
 c_k =\frac{\la_k}{\pi_k(\la_k)}\,,
\ee
 where $\pi_k(\la)$ is defined by
\bel{pikdef}
 \pi_k(\la) = \sum_{j\ge k-1}\frac{e^{-\la}\la^j}{j!}\,,
\ee
 and $\la_k$ is the positive number minimising the right hand side
of~\eqn{ckdef}. Recently, a number of simpler proofs establishing
the threshold $c_k/n$ for the $k$-core have been published (see
Kim~\cite{K}, Cain and Wormald~\cite{CW}, and Janson and
Luczak~\cite{JL}). \bs

In what follows, we write a.a.s.\ to denote an event which occurs
with probability tending to one as $n \to \infty$. In comparison to
studying the $k$-core in random graphs, it appears to be
substantially more difficult to analyse the appearance of
$k$-regular subgraphs when $k \geq 3$. One reason is that it is
NP-hard to determine whether a graph contains such a subgraph, and
there is no analogue of the simple vertex deletion algorithm which
produces the $k$-core. As every $k$-regular subgraph is contained in
the $k$-core, we deduce that $\G(n,p)$ a.a.s.\ does not contain a
$k$-regular subgraph whenever $p$ is below the threshold for the
$k$-core described in (\ref{ckdef}) and (\ref{pikdef}).
Bollob\'{a}s, Kim and Verstra\"{e}te~\cite{BKV} showed that
$\G(n,p)$ a.a.s.\ contains a $k$-regular subgraph when $p$ is,
roughly, larger than $4c_k/n$, and conjectured a sharp threshold for
the appearance of $k$-regular subgraphs in $\G(n,p)$. In the same
paper it was shown that for some $c > c_3$, the $3$-core of
$\G(n,c/n)$ has no 3-regular subgraph a.a.s., whereas for $c \geq
c_4$, the $4$-core of $\G(n,c/n)$ contains a 3-regular subgraph
a.a.s. In support of the conjecture of a sharp threshold, Pretti and
Weigt~\cite{PW} numerically analysed equations arising from the
cavity method of statistical physics to conclude empirically that
indeed, there is a sharp threshold for the appearance of a
$k$-regular subgraph of a random graph. For $k > 3$ they concluded
that it is the same as the threshold for the $k$-core, which is at
odds with~\cite[Conjecture 1.3]{BKV}. For $k = 3$, these thresholds
differ, as shown using the first moment method in~\cite{BKV}.

\bs

In this paper, we improve the window of the threshold for
$k$-regular subgraphs in $\G(n,p)$ by proving Theorem~\ref{main}
below. A {\it $k$-factor} of a graph is a spanning $k$-regular
subgraph, and a graph is {\it $k$-factor critical} if, whenever we
delete a vertex from the graph, we obtain a graph which has a
$k$-factor.

\begin{theorem}\lab{main}
There exists an absolute constant $k_0$ such that for $k \geq k_0$,
the $(k + 2)$-core of a random graph $\G(n,p)$ is nonempty and
contains a $k$-factor or is $k$-factor-critical a.a.s.
\end{theorem}

Theorem \ref{main} will be proved in Section \ref{thm:main}. We
remark that the first nonempty $k$-core of the random graph process
a.a.s.\ contains many vertices of degree $k + 1$ adjacent to $k+1$
vertices of degree $k$, so the $k$-core cannot contain a $k$-factor
and cannot be $k$-factor critical a.a.s. Bollob\'{a}s, Cooper,
Fenner and Frieze~\cite{BCFF} conjectured that the $(k + 1)$-core
contains $\lfloor k/2 \rfloor$ edge disjoint hamiltonian cycles
a.a.s., so Theorem~\ref{main} supports this conjecture.

\bs

The value of $c_k$ can be determined approximately for large $k$ as
follows. This corrects, and sharpens, the error term of the formula
given in~\cite{PSW}. All logarithms in this paper are natural, and
$\mathbb N$ is the set of positive integers.
\begin{lemma}\lab{ckasy}
For any $k \in \mathbb N$, let $q_k =  \log k - \log (2\pi )$. Then
$$
c_k=k+(k q_k)^{1/2}+\Bigl(\frac{k}{q_k}\Bigr)^{1/2}
+ \frac{q_k -1}{3} + O\Bigl(\frac{1}{\log k}\Bigr) \quad \mbox{\rm
as } k \rightarrow \infty.
$$
\end{lemma}

Lemma~\ref{ckasy} is proved in Section~\ref{lemproof}. It follows
immediately from this lemma that
$$c_{k+2} = c_k+2+O\Bigl(\frac{1}{\log k}\Bigr) \,. $$
Hence, combining the lemma and Theorem~\ref{main}, we have pinned
down the threshold for the appearance of $k$-regular subgraphs in
$\G(n,p)$ to a window for $p$ of width $2/n + O(1/n\log k)$. The
following questions remain: (1) to determine whether there is a
sharp threshold for the appearance of a $k$-regular subgraph, and
(2) whether the $(k+1)$-core of a random graph, when it is a.a.s.\
nonempty, contains a $k$-factor or is $k$-factor critical a.a.s. \bs

Throughout the paper, we denote by $\G(n,p)$ the Erd\H{o}s-R\'{e}nyi
model of random graphs. If $G$ is a graph with vertex set $V(G)$,
then $\lambda(S,T)$ denotes the number of edges of $G$ with one
endpoint in $S$ and one endpoint in $T$, where $S,T \subseteq V(G)$.
If $S = T$, we write $\lambda(S)$ instead. The number of components
of a graph $G$ is denoted by $\omega(G)$.

\section{Factors of Graphs}

In this section, we allow graphs to contain multiple edges. Let $G$
be a graph and let $k \in \mathbb N$. Recall that a {\it $k$-factor}
of $G$ is a spanning subgraph of $G$ all of whose vertices have
degree $k$. A graph is {\it $k$-factor-critical} if the deletion of
any vertex of $G$ results in a graph with a $k$-factor. In
particular, a {\it 1-factor} of $G$ is a perfect matching of $G$.
Tutte's $1$-Factor Theorem gives the following necessary and
sufficient condition for a graph $G$ to have a 1-factor.

\begin{theorem} \lab{thmtutte}
Let $G$ be a graph, and let $\mbox{\rm o}(G)$ denote the number of
components of $G$ with an odd number of vertices. Then $G$ has a
1-factor if and only if \be \label{tutte} \mbox{\rm o}(G - X) \leq
|X| \ee for every set $X \subseteq V(G)$.
\end{theorem}

Using Tutte's 1-Factor Theorem applied to a new graph, $\phi(G)$, a
necessary and sufficient condition can be found for a graph to have
a $k$-factor (this is a special case of Tutte's $f$-factor theorem;
see Lov\'{a}sz and Plummer~\cite{LP} for details). To construct
$\phi(G)$, let $V(G) = \{v_1,v_2,\dots,v_n\}$, and let $V(\phi(G)) =
U \cup V$, where $U$ and $V$ are disjoint sets and $V$ is
partitioned into independent sets $(V_1,V_2,\dots,V_n)$ with $|V_i|
= d(v_i)$ and $U$ is partitioned into sets $(U_1,U_2,\dots,U_n)$
such that $|U_i| = k$. Then $\phi(G)$ consists of all edges between
$U_i$ and $V_i$ for $i \in \{1,2,\dots,n\}$, together with a
matching on $V$ such that when we contract all the independent sets
$V_i$ in $\phi(G) - U$ to single vertices, we obtain $G$. An example
is shown below in Figure 1, where $G$ is a quadrilateral, $k = 2$,
and a $1$-factor in $\phi(G)$ corresponds to a $2$-factor in $G$.

\bs

\SetLabels
\endSetLabels
\centerline{\AffixLabels{\includegraphics[width=3in]{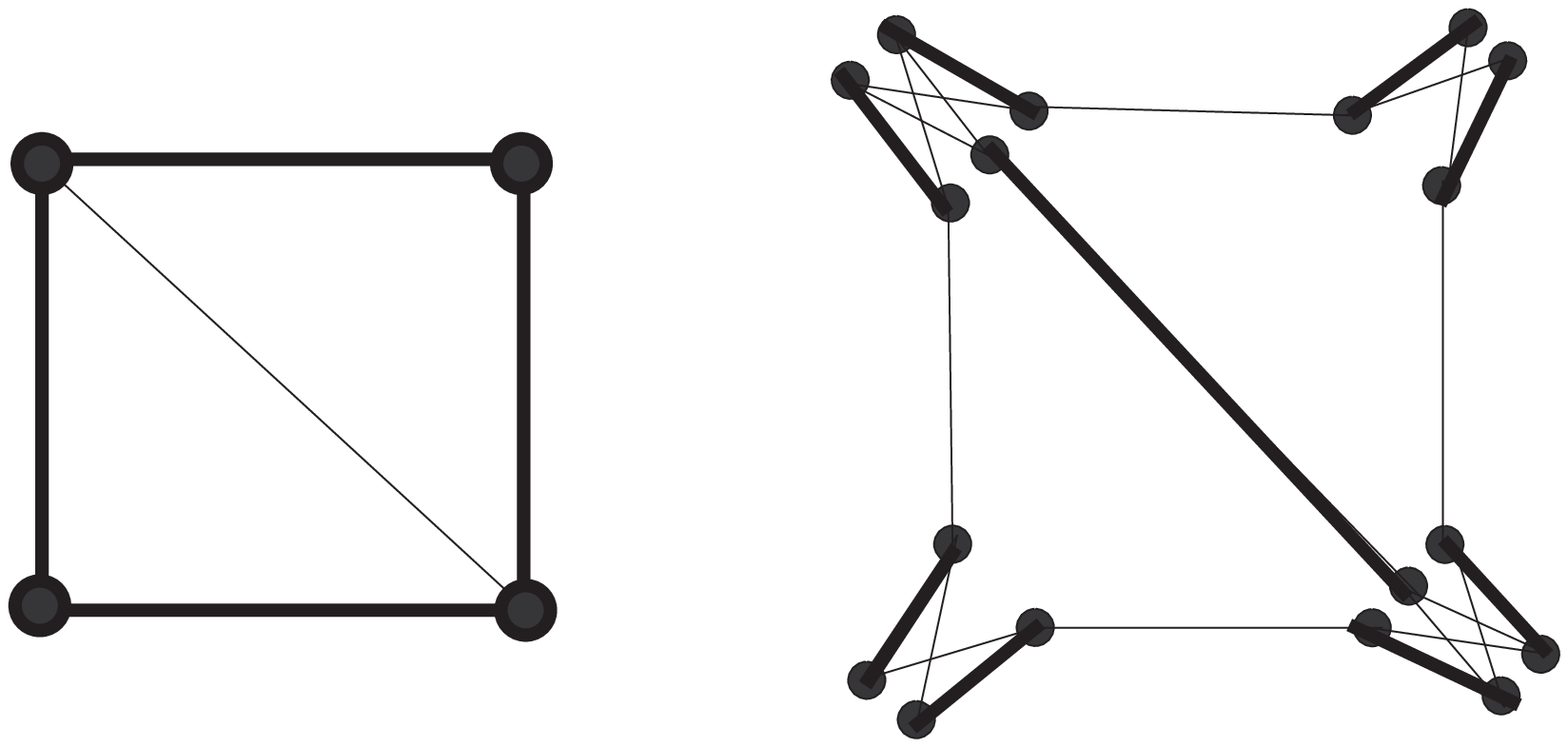}}}

\begin{center}
{\sf Figure 1 : 2-factor in $G$ and corresponding 1-factor in
$\phi(G)$.}
\end{center}

\bs

It is not hard to see that $G$ has a $k$-factor if and only if
$\phi(G)$ has a $1$-factor. We will show that the following
condition on $G$ is enough to guarantee a $k$-factor in $G$. It is
convenient to define
\[ \delta_k(G) = \left\{\begin{array}{ll}
0 & \mbox{ if }k|V(G)|\mbox{ is even} \\
1 & \mbox{ if }k|V(G)|\mbox{ is odd}.
\end{array}\right.\]

\begin{lemma}\lab{factor} Let
$k \in \mathbb N$, and let $G$ be a connected graph such that for
every pair of disjoint sets $S,T \subseteq V(G)$ for which $S \cup T
\neq \emptyset$,
\begin{equation}\lab{kfactor}
\sum_{v \in T} d(v) + k|S| \geq  \omega\big(G - (S \cup T)\big) +
k|T| + \lambda(S,T) + \delta_k(G)\,.
\end{equation}
Then $G$ has a $k$-factor or is $k$-factor critical according as
$\delta_k(G) = 0$ or $\delta_k(G) = 1$.
\end{lemma}
\proof We first consider the case $\delta_k(G) = 0$. Let $H =
\phi(G)$. To show that $G$ has a $k$-factor, it is sufficient to
show that for all $X \subseteq V(H)$, $\mbox{o}(H - X) \leq |X|$, by
\eqn{tutte}. If $X = \emptyset$, then this follows from the fact
that $H$ is connected and
$$|V(H)| = k|V(G)| + \sum_{i=1}^n d(v_i) \equiv k|V(G)| \equiv 0
{\rm \ mod\ } 2.$$
 In what remains, we verify \eqn{tutte} for $X \neq \emptyset$.
Suppose that for some $i$, $0< | X \cap U_i |<|U_i|$. Since $U_i$
and $V_i$ form a complete bipartite graph, we may delete one of the
vertices of $U_i$ from $X$, and the number of components of $H - X$
does not change. Then the right side of~\eqn{tutte} decreases, and
the left side decreases by at most 1. Hence we assume $X \cap U_i =
\emptyset$ or $X \cap U_i = U_i$ for all $i$.

\bs

Define the following sets of vertices of $G$:
\begin{eqnarray*}
S = \{v_i \in V(G) : X \cap U_i = U_i\}, \; \; \; \; T = \{v_i \in
V(G) : X \cap V_i = V_i\},
\end{eqnarray*}
and the following sets of vertices of $H$:
\begin{eqnarray*}
  Y = \bigcup_{i=1}^n \; \{V_i - X : v_i \in S\}, \; \; \;   X_0 = \bigcup_{i=1}^n \;
\{X \cap V_i : v_i \not \in T\}, \; \; \;  X_1 =  \bigcup_{i=1}^n \;
\{V_i : v_i \in T\}.
\end{eqnarray*}
Suppose that for some $i$, $v_i \in S \cap T$ (that is, $X\cap
U_i=U_i$ and $X\cap V_i=V_i$). But then we may delete all the
vertices of $U_i$ from $X$, and both the left side and the right
side of~\eqn{tutte} decreases by $k$. Hence we assume $S \cap T =
\emptyset$.

For convenience, given a subgraph $F$ of $\phi(G)$, we write
$\phi^{-1}(F)$ for the subgraph of $G$ obtained from $F$ by deleting
all vertices of $U$ and contracting all the sets $V_i\cap V(F)$. For
each component $F$ of $H - X$ containing a vertex of $U = \bigcup_{i
\leq n} U_i$, either $F$ is an isolated vertex in $U$, or
$\phi^{-1}(F) - S$ is a component of $G - (S \cup T)$. So at most
$k|T| + \omega\big(G - (S \cup T)\big)$ components of $H - X$
contain a vertex of $U$. Now let $F$ be a component of $H - X$
containing no vertices of $U$. Then $|V(F)| \leq 2$, so the only
components of $H - X$ containing no vertices of $U$ and contributing
to $\mbox{o}(H - X)$ are isolated vertices -- and these are vertices
of $Y$ which are adjacent to a vertex of $X = X_0 \cup X_1$. So the
number of these isolated vertices is:
\[ \lambda_H(Y,X_0 \cup X_1) = \lambda_H(Y,X_0) + \lambda_H(Y,X_1) = \lambda_H(Y,X_0) + \lambda(S,T) \leq
|X_0| + \lambda(S,T)\] where $\lambda_H(Y,X)$ is the number of edges
between $X$ and $Y$ in $H$. It now follows using~(\ref{kfactor})
that
\begin{eqnarray*}
\mbox{\rm o}(H - X) &\leq& \omega(G - (S \cup T)) + k|T| + \lambda(S,T) + |X_0| \\
                       &\leq& \sum_{v \in T} d(v) + k|S| + |X_0| \; \; = \; \; |X_1| + k|S| + |X_0| \; \;
= \; \;  |X|.
\end{eqnarray*}
So we have shown $|X| \geq \mbox{o}(H - X)$ for every set $X
\subseteq V(H)$, as required. This completes the proof for
$\delta_k(G) = 0$.

\bs

Finally, suppose $\delta_k(G) = 1$ and let $v \in V(G)$. We show
that $G - \{v\}$ has a $k$-factor. Let $G'$ be the graph obtained by
adding a vertex $v'$ to $G$ and joining $v$ to $v'$ with $k$
parallel edges. Then $G'$ is connected, and $\delta_k(G') = 0$.
Furthermore, it is straightforward to check that \eqn{kfactor} is
satisfied in $G'$. By the first part of the proof, $G'$ has a
$k$-factor. Deleting both $v$ and $v'$ from this $k$-factor of $G'$,
we get a $k$-factor in $G - \{v\}$, as required. \qed

\bs

\section{Structure of the $k$-core}

In this section we describe the structure of the $k$-core in
$\G(n,p)$; this material will be used throughout the proof of
Theorem \ref{main}. We will assume throughout that $p = c/n$ where
$c > c_k$, so that the $k$-core of $\G(n,p)$ is a.a.s.\ nonempty. We
let $K$ denote this nonempty $k$-core.

\bs

In the first lemma, $\partial X$ denotes the set of edges of $K$
with exactly one endpoint in a set $X \subset V(K)$. The lemma seems
to be well known, and follows, for example, from Benjamini, Kozma
and Wormald~\cite[Lemma 5.3]{bkw}. (That lemma concerns graphs with
a given degree sequence, all degrees between 3 and $n^{0.02}$. See
the proof of~\cite[Theorem 4.2]{bkw} to find the connection with the
following.)

\begin{lemma}\lab{expansion}
There is a positive constant $\gamma$ such that the following holds.
Fix $k \geq 3$. Then a.a.s.\ every set $X \subset V(K)$ of at most
$\frac{1}{2}|V(K)|$ vertices satisfies
\[ |\partial X| \geq \gamma k |X|.\]
\end{lemma}

Throughout the rest of the paper, $\gamma$ denotes the constant
appearing in Lemma~\ref{expansion}.

\begin{lemma}\lab{kcorestructure}
Let $k > 2/\gamma$. Then a.a.s.\ for every set $Y \subset V(K)$ of
size at most $s(n) = \log n/2ec\log\log n$, $K - Y$ contains a
component with more than $|V(K)| - 2s(n)$ vertices.
\end{lemma}
\begin{proof}
We first show that a.a.s.\ there are no sets of $2y$, $y \le s(n)$
 vertices in $\G(n,p)$ inducing at least $2y + 1$ edges: the
expected number of subgraphs of $2y$ vertices with at least $2y + 1$
edges, for some $y \leq s(n)$, is at most
\begin{eqnarray*}
\sum_{y \leq s(n)} {n \choose 2y}{{2y \choose 2} \choose 2y + 1}
p^{2y + 1}  &<& \sum_{y \leq s(n)} \frac{(\log
n)^{2y+1}}{n} \\
&<& n^{-1/2}.
\end{eqnarray*}
So the claim follows from Markov's inequality. Thus, we may  assume
that   all sets of $2y \leq 2s(n)$ vertices in $K$ induce at most
$2y$ edges. We may also assume that a.a.s.\ the property in
Lemma~\ref{expansion} holds.

It suffices to show that if $X$ is the vertex set of a union of
components of $K - Y$ and $|X| \leq \frac{1}{2}|V(K)|$, then $|X| <
s(n)$. Suppose $|X| \geq s(n)$. From the property in Lemma~\ref{expansion},
\[ |\partial X| \geq \gamma k|X|. \]
Suppose $|Y| = y$. By averaging, some $Z \subset X$ of size $|Y|
\leq s(n)$ satisfies
\[ \lambda(Y \cup Z) \geq \gamma ky > 2y.\]
However $|Y \cup Z| = 2|Y| = 2y$, which is a contradiction.
\end{proof}

\bs

In fact \L uczak~\cite{L} showed that the $k$-core is $k$-connected
a.a.s., as stated in the next lemma.

\begin{lemma}\lab{connected}
For $k \geq 3$ and $c > c_k$, the $k$-core of $\G(n,c/n)$ is
$k$-connected a.a.s.
\end{lemma}

Our final lemma is a large deviation result for the degrees of the
vertices of the $k$-core. Essentially, the degree of a vertex in $K$
has (asymptotically) a truncated Poisson distribution, which gives a
precise bound on the number of vertices which deviate from degree
$c$ in $K$.

\begin{lemma}\lab{new1}
For all $\eps>0$ there exists $k_\eps$ such that for $k>k_\eps$  and
$c_k<c< 2k$, it is a.a.s.\ true that
$$
|d(v)-c|\ge \eps \sqrt{k\log k}
$$
for at most $ \eps |K|$ vertices $v$ of $K$.
\end{lemma}
\proof Let $\eps > 0$, and fix $k$, and $j\ge k+2$.
From~\cite[Corollary 3 and Erratum]{CW}, if $n_j$ denotes the number
of vertices of degree $j$ in $\K$, then a.a.s.\
 \bel{njasy}
n_j = \frac{e^{-\mu}\mu^j}{j!}n +o(n),
 \ee
where $\mu = \mu_{k,c}$ is the larger of the two positive solutions
of the equation
 \bel{muceq}
\frac{\mu}{c}= e^{-\mu} \sum_{i\ge k-1}\frac{\mu^i}{i!}.
 \ee
(The fact that there are two such solutions is known to be
guaranteed by the fact that $c>c_k$.) Let $\eps_1>0$, and suppose
that $\mu=\Theta( k)$. Then, since the Poisson distribution is
asymptotically normal with variance equal to its mean, we have for
sufficiently large $k$
 \bel{poisson}
\sum_{|i-\mu|\ge \eps_1\sqrt{k\log k} }e^{-\mu}  \frac{\mu^i}{i!} <
\eps_1.
 \ee
Also, by Lemma~\ref{ckasy}, we may assume that $c>k+\frac12\sqrt{k
\log k}$. Suppose that $c-2$ is substituted for $\mu$ in~\eqn{muceq}. It is then elementary to obtain
that the right hand side of~\eqn{muceq} is greater than $1-
1/k$. Recalling also that $c<2k$, this is greater than the left hand
side of~\eqn{muceq}. On the other hand, if anthing larger than $c$ is substituted for $\mu$ in~\eqn{muceq} then the left hand side is greater than the right, since the right is equal to a probability strictly less then 1.  So by continuity, $c-2 < \mu_{k,c}<c$. Taking
$\eps_1$ slightly smaller than $\eps$, the lemma follows
from~\eqn{poisson} and~\eqn{njasy}. \qed

\section{Proof of Theorem \ref{main}}\label{thm:main}

In this section, we denote by $K$ the $(k+2)$-core of $\G(n,p)$
where $pn = c$ and $c_{k+2} < c$, for $k \geq 3$. To prove Theorem
\ref{main}, we show that there exists $k_0$ such that for $k \geq
k_0$, (\ref{kfactor}) holds in $K$ a.a.s. The value of $k_0$ will
not be optimized in the proof to follow. To prove the theorem, we
consider a number of cases according to the sizes of the sets $S$
and $T$ in the lemma, where $S \cup T \neq \emptyset$. It is
convenient throughout to let $s(n) = \log n/2ec\log \log n$.

\ms

\ms \no {\bf Case 1\ } {\it $\displaystyle{|S| + |T| < s(n)}$.}

\ms

Let $Y = S \cup T$. By Lemma \ref{kcorestructure}, $K - Y$ contains
a component with more than $|V(K)| - 2s(n)$ vertices a.a.s. Let
$\omega_{\rm s}(K - Y)$ denote the number of components of $K - Y$
of size less than $s(n)$, and let $X$ be the set of vertices in
these components. As in the proof of Lemma~\ref{kcorestructure},
$\lambda(X) \leq |X|$ and $\lambda(X \cup Y) \leq |X| + |Y|$ a.a.s.
However, every vertex of $X$ has degree at least $k + 2$, by
definition of $K$, so $\lambda(X) + \lambda(X,Y) \geq (k+1)|X|$
a.a.s.\ since $\lambda(X) \leq |X|$ holds a.a.s. It follows that
\begin{eqnarray*}
|X| + |Y| \; \; \geq \; \;  \lambda(X \cup Y) &=& \lambda(X) + \lambda(X,Y) + \lambda(Y) \\
&\geq& k|X| + |X| + \lambda(S,T)\\
&\geq& k \cdot \omega_{\rm s}(K - Y) + |X| + \lambda(S,T).
\end{eqnarray*}
Therefore, since $\omega_{\rm s}(K - Y) = \omega(K - Y) - 1$,
\begin{eqnarray*}
k \cdot \omega\big(K - (S \cup T)\big) - k + \lambda(S,T) + k|T|
&\leq& |S| + (k+2)|T| - |T| \,.
\end{eqnarray*}

Since $k \geq 3$, and $\omega\big(G - (S \cup T)\big) \geq 1$,
\begin{eqnarray*}
\omega\big(K - (S \cup T)\big) + \lambda(S,T) + k|T|
&\leq& k \cdot \omega\big(K - (S \cup T)\big) - k+1 + \lambda(S,T) + k|T|\\ \\
&\leq& |S| + (k+2)|T| - |T| \\ \\
&\leq& k|S| + \sum_{v \in T} d(v) - (k - 1)|S| - |T| +1\\ \\
&\leq& k|S| + \sum_{v \in T} d(v) - 2|S| - |T| + 1 \,.
\end{eqnarray*}
If $S \neq \emptyset$ or $|T| > 1$, then this is less than $k|S| +
\sum_{v \in T} d(v)$, as required. If $S=\emptyset$ and $|T|=1$,
then $\omega(K - T) = 1$ a.a.s.\ (see Lemma \ref{connected}), which
implies
\[ \omega(K - T) + k|T| = 1 + k < 2 + k \leq \sum_{v \in T} d(v).\]
Therefore (\ref{kfactor}) holds a.a.s.\ in $K$, so $K$ has a
$k$-factor or is $k$-factor critical, by Lemma \ref{factor}.

\bs

For the rest of the proof, $\eps_0$ is an absolute constant; we will
take $\eps_0 = e^{-9}$ for definiteness.

\bs

\no {\bf Case 2\ } $|S| + |T| \geq s(n)$, $|T| < \eps_0 n$, and $|S|
< 4\eps_0 n$.

\bs

Let $Y = S \cup T$. In this case we estimate $\omega(K - Y)$ and
$\lambda(S,T)$ separately. First we show that $\omega(K - Y) \leq
|Y|/2$ a.a.s.\ provided that $k$ is large enough to ensure
that $\gamma(k+2) \geq 4$ (so this tells us $k_0 \geq 4/(\gamma -
2)$ is required in our proof). It suffices to show that if $X$ is
the vertex set of any union of components of $K - Y$, then a.a.s.\
$|X| < |Y|/2$ or $|X| > n/2$. Suppose that $|Y|/2 \leq |X| \leq
n/2$. Then Lemma~\ref{expansion} shows $\lambda(X,Y) \geq \gamma(k +
2)|X|$. Let \[ I = \{y : s(n) \leq y \leq 5\eps_0 n\} \; \; \mbox{
and } \; \; I_y = \Big\{x : \frac{y}{2} \leq x \leq
\frac{n}{2}\Big\}.\] The expected number of pairs of sets $(X,Y)$ in
$\G(n,p)$ satisfying the above requirements is
\begin{eqnarray*}
\sum_{I_y \times I} {n \choose x+y} {x + y \choose x} {xy \choose
\gamma(k+2)x} p^{\gamma (k + 2)x} &<& \sum_{I_y \times I} e^{3x +
\gamma(k+2)x}2^{x+y} \Bigl(\frac{n}{x}\Bigr)^{3x}
\Bigl(\frac{y}{n}\Bigr)^{\gamma(k+2)x} \\ \\
&<& \sum_{I_y \times I} e^{3\gamma(k + 2)x/2}2^{3x}
\Bigl(\frac{y^2}{xn}\Bigr)^{ \gamma(k+2)x/2} \\ \\
&<& \sum_{I_y \times I} e^{-\gamma(k+2)x} e^{3x} \; \quad (\mbox{since }y \leq 5\eps_0 n)\\ \\
&\leq& \sum_{I_y \times I} e^{-x} \quad \quad \quad \quad \quad (\mbox{since }\gamma(k+2) \geq 4)\\ \\
&<& \int_{I_y \times I} e^{-(x-1)} dx dy \\ \\
&<& \int_{I} \; \; \; e^{- (y-2)/2} dy \\ \\
&<& e^{- (s(n) - 2)/2} \; \; = \; \;  o(1) \,,
\end{eqnarray*}
where the sums are over $(x,y) \in I_y \times I$. So in fact the
expected number of sets $X$ and $Y$ as described above is $o(1)$. By
Markov's Inequality, we conclude that $\omega(K - Y) \leq |Y|/2$
a.a.s.

\ms

To finish verifying inequality~\eqn{kfactor} in case 2, it remains
to show that a.a.s.
\begin{equation}\label{stbound}
\lambda(S,T) < \frac{3}{2}|T| + \Bigl(k- \frac{1}{2}\Bigr)|S|.
\end{equation}
Let $|S| = \sigma n$ and $|T| = \tau n$ and let $\rho = (k -
\frac{1}{2})\sigma + \frac{3}{2}\tau$ (here $\sigma$ and $\tau$ are
allowed to depend on $n$). Then the number of ways of choosing the
sets $S$ and $T$ is bounded above by
\[ {n \choose \sigma n}{n \choose \tau n} < \Bigl(\frac{\sigma}{e}\Bigr)^{-\sigma n}
\Bigl(\frac{e}{\tau}\Bigr)^{-\tau n}.\] The probability that there
are at least $\rho n$ edges between $S$ and $T$ is at most
 \bean
{\sigma \tau n^2 \choose \rho n}\left(\frac{c}{n}\right)^{\rho n }
&\le& \left(\frac{e\sigma \tau c}{\rho}\right)^{\rho n}\\
&=& \left(\frac{e c}{ (k-1/2)/\tau +
3/2\sigma }\right)^{(k- 1/2)\sigma n +
3\tau n/2}.
 \eean
Multiplying by the bound on the number of choices of $S$ and $T$,
the bound (\ref{stbound}) is true a.a.s.\ if
$$
\left(\frac{e c}{ (k-1/2)/\tau +
3/2\sigma }\right)^{(k- 1/2)\sigma n +
3\tau n/2} < \left(\frac{\sigma}{e}
\right)^{\sigma}\left(\frac{\tau}{e} \right)^{\tau}
$$
for which it suffices that
\begin{equation}\label{suffices}
 \frac{e c}{ (k-1/2)/\tau +
3/2\sigma }  <\min \Bigl\{ \Bigl(\frac{\sigma}{e}\Bigr)^{1/(k-
1/2)}, \Bigl(\frac{\tau}{e}\Bigr)^{2/3} \Bigr\}.
\end{equation}
Since $\tau < e^{-9}$ and $c < 2k - 1$ for large enough $k$,
\[
\frac{e c}{ (k-1/2)/\tau + 3/2\sigma } < \frac{e\tau c}{k-1/2} <
2e\tau < \Bigl(\frac{\tau}{e}\Bigr)^{2/3}.\] To prove that the
expression on the left in (\ref{suffices}) is less than
$(\sigma/e)^{1/(k-1/2)}$,   we can assume
\[ \Bigl(\frac{\sigma}{e}\Bigr)^{1/(k-1/2)} < \Bigl(\frac{e\tau c}{k-1/2}\Bigr) < e^{-7}.\]
It follows that $\sigma < e^{-4k}$ and therefore, since $2k - 1 \geq
2$,
\[
\frac{e c}{ (k-1/2)/\tau + 3/2\sigma } < e\sigma c < e^{-2k}c
\sigma^{1/2} < \sigma^{1/2} \leq \sigma^{1/(k-1/2)}.\]

\bs

\no {\bf Case 3\ } $|S| + |T| \geq s(n)$ and $|T| < \eps_0 n$ and
$|S| \ge 4\eps_0 n$.

\bs

To prove that the requirements of Lemma~\ref{factor} are satisfied
a.a.s.\ for $k$ sufficiently large, it is enough to show that
a.a.s.\ for every pair of sets $(S,T)$ under consideration,
$\lambda(S,T) \le \frac34 k|S|$. This is because the sum of vertex
degrees in $T$ is at least $k|T|$, and also $\omega\big(G - (S \cup
T)\big)\le n< \frac14 k|S|$ for large enough $k$ (as $\eps_0$ is
fixed). We may assume that $c\sim c_k$, and so if $k$ is
sufficiently large, Lemma~\ref{ckasy} shows that $c<\frac32 k$. Thus
it suffices to show $\lambda(S,T) \le   2c \eps_0 n$, since this is
at most $ \frac12 c|S|$. This follows immediately once we show that
a.a.s.\ all sets of at most $\eps_0 n$ vertices (in particular, $T$)
have total degree at most $2c \eps_0 n$.

It is well known that in the random graph $\G(n,c/n)$, the vertex
degrees are have  asymptotically Poisson distribution with mean $c$:
the number of vertices of degree $j$ is a.a.s.\ asymptotic to
$e^{-c}c^j/j!$. It follows that the sum of the degrees of those
vertices of degree less than $3c/2$ is a.a.s.\ asymptotic to
$$
n\sum_{j<3c/2}j\frac{e^{-c}c^j}{j!}.
$$
Since the Poisson distribution is asymptotically normal with mean
$c$, for large enough $c$ (i.e.\ large enough $k$) we have
$$
\sum_{j<3c/2}j\frac{e^{-c}c^j}{j!} > c-\eps_0.
$$
Since the sum of vertex degrees is a.a.s.\ asymptotic to $cn$, the
ones of degree at least $3c/2$ a.a.s.\ have total degree less than
$\eps_0n$. Assuming this is true, there are at most $2\eps_0n/3c$
such vertices, and any set of at most $\eps_0n$ vertices thus has
total degree at most
$$
 \eps_0n  +  \frac{3c}{2}\eps_0n <2c\eps_0n
$$
as required.

\bs

\no {\bf Case 4\ } $ |T| \ge \eps_0 n $.

\bs



By Lemmas~\ref{ckasy} and~\ref{new1}, for all $\eps>0$, if $k$ is
sufficiently large, then a.a.s.
$$
\sum_{v \in T} d(v) > \big(k+ (1-\eps) \sqrt{k \log k} \big) |T|.
$$
So by Lemma~\ref{factor} and the fact that $n=O(|T|)$, it is enough
to show, for some $\eps>0$, that a.a.s.
 \bel{enough2}
(1-\eps) \sqrt {k \log k}\, |T| + k|S| \ge \lambda(S,T).
 \ee
We will prove this by considering the cases $|S| < \eta n$ and $|S|
\geq \eta n$ separately, where $\eta = \frac{1}{4}\eps_0$. For $|S|
< \eta n$, we will use
$$\lambda(S,T)\le  \sum_{v\in S}d(v). $$
For this, we may assume that if $|S|= \sigma n$ then $S$ contains
the $\sigma n$ vertices of largest degree in $G$ (and that they have
the same degrees in $K$). Using the argument about the degrees of
$G\in\G(n,p)$ as in Case~3, it is straightforward to show that
a.a.s.\ these vertices have total degree at most
 \bean
c\sigma n+O(\sqrt c\,  n)&  < & \big(k + 2  \sqrt{k \log k}\big)|S| +O(\sqrt k\,  n)\\
&< &k |S|+ 3 \sqrt{k \log k}\,  |S|
 \eean
for $k$  sufficiently large, which is at most
$$
  (1-\eps) \sqrt {k \log k}\, |T| + k|S|
$$
since $|S|< \eta n \leq \frac14|T|$. This gives~\eqn{enough2}.

It only remains to treat those sets $S$ for which $|S| \ge \eta n$.
Then using the same argument as with $T$ in Case~3, the sum of degrees
of vertices in $S$ is a.a.s.\ at most $\big(k+(1+\eta)\sqrt{k \log
k}\big)|S|$.

For a set $S$ of this size in $\G(n,p)$, the expected value of
$\lambda(S)$ is
$$
{\eta n \choose 2} \frac {c}{n} \sim \frac 12   c n \eta^2 > \frac12
kn\eta^2
$$
since $c > k$. Moreover, $\lambda(S)$ is binomially distributed. So
by Chernoff's inequality (see for example~\cite[Theorem 2.1]{JLR}),
$$
\pr\Bigl(\lambda(S)\le \frac{kn\eta^2}{4}\Bigr) \le
e^{- kn\eta^2/16 } =o(2^{-n})
$$
for sufficiently large $k$ (recall that $\eta$ is an absolute
constant). Hence, a.a.s.\ every set $S$ that is this large induces a
subgraph of at least $\frac{1}{4}kn\eta^2 \ge \frac{1}{4}|S|
k\eta^2$ edges. Provided $S\subseteq V(K)$, it contains exactly the
same number of edges in $K$ as in $\G(n,p)$. Hence we have that
a.a.s.\ for all such $S$ and $T$,
 \bean
\lambda(S,T) &\le& \sum_{v\in S} d(v) - \frac{1}{2}|S|  k\eta^2 \\
\\
&\le &  \big(k+(1+\eta)\sqrt{k  \log k}\big)|S| - \frac{1}{2}|S|
k\eta^2\\ \\
& \le& k|S|
 \eean
for large enough $k$. This gives~\eqn{enough2}, as required. \qed

\section{Proof of Lemma~\ref{ckasy}}\label{lemproof}

A weakened version of the main result in Pittel, Spencer and
Wormald~\cite{PSW} is that if $c$ is fixed, $\G(n,c/n)$ a.a.s.\ has
no $k$-core if  $c< c_k $, a.a.s.\ has one if $c> c_k$, where $c_k$
is defined in (\ref{ckdef}). A little calculation shows that $c_k$
and $\la_k$ ($\la_k$ is also defined in  (\ref{ckdef})) satisfy
 \bea
 c_k \pi_k(\la_k)& =& \la_k  \lab{clam1}\\
 c_k&=&(k-2)!e^{\la_k}\la_k^{-(k-2)} \lab{clam2}
 \eea
with $\pi_k$ defined in~\eqn{pikdef}. Substituting~\eqn{clam2}
and~\eqn{pikdef} into~\eqn{clam1} gives
$$
\la_k = \sum_{j\ge 0} \frac{{\la_k}^{j+1}}{[k+j-1]_{j+1}}
$$
(where square brackets denote falling factorials) and so,
multiplying by $(k-1)/\la_k$, we obtain
$$
k-1 = \sum_{j\ge 1} \frac{{\la_k}^{j}}{[k+j-1]_{j}}.
$$
Since the right hand side is an increasing function of $\la_k$, the
value of $\la_k$ is uniquely determined. Moreover, since $(k +
j)/\lambda_k$ is exactly the ratio of the $j$th to the $(j + 1)$th
term in the summation, the largest term
in the summation occurs for $k+j \approx \la_k$  and from elementary considerations it is easy to see that $\la_k=
k+O(\sqrt{k\log k})$. Thus, putting \bel{tdef} \la_k=(k-2)(1+t), \ee
we know that $t=o(1)$. In addition, rewriting~\eqn{pikdef} as
\bel{pikbackwards}
 \pi_k =
1-\sum_{j\le k-2}\frac{e^{-\la}\la^j}{j!}, \ee we now see that
$\pi_k=1-o(1)$, and hence also $c_k\sim k$.

\bs

To get a slightly better bound on $t$ straight away,
substitute~\eqn{clam2} into~\eqn{clam1}, use Stirling's formula with
its correction term due to Robbins: $j!=(j/e)^j\sqrt{2\pi j}
(1+O(1/j))$, and take logarithms to give \bel{logpik} \log \pi_k =
\frac12 \log\Bigl(\frac{k-2}{2\pi}\Bigr)+ (k-1)\log(1+t) -(k-2)t
+O\Bigl(\frac{1}{k}\Bigr). \ee Recalling from above that $\log
\pi_k=o(1)$ and $t=o(1)$, we may expand $\log(1+t)$ to show that
 \bel{tasy}
t\sim \Bigl(\frac{q_k}{k}\Bigr)^{1/2}
 \ee
where $q_k = \log k - \log(2\pi)$.

Taking out a factor of $1/c_k$ from the terms in the summation
in~\eqn{pikbackwards}, using~\eqn{clam2} we obtain
$$
\pi_k=1-\frac{1}{c_k}
\sum_{m=0}^{k-2}(1+t)^{-m}\left(\frac{k-1}{k-2}\right)^m\prod_{j=1}^m\left(1-\frac{j}{k-2}\right).
$$
The terms in the summation are monotonically decreasing. Since
$(1+t)^{-m} =\displaystyle{\exp(-mt+O(mt^2))}$, we see that, for any
$\eps>0$, the terms for $m>k^{1/2+\eps}$ sum to $o(1/k)$. For
$m=O(k^{1/2+\eps})$, we see after expanding that the product
over $j$ is
$$
 e^{- m^2/2k  +O( m/k + m^3/k^2 )} =  1+ O\Bigl(\frac{m^2}{k}  + \frac{m^3}{k^2}\Bigr).
$$
Putting $r=\log(1+t)$ and recalling $c\sim k$, we now have
$$
\pi_k \; \; = \; \;  1 \; \; - \; \; \frac{1}{c_k}
\sum_{m=0}^{k^{1/2+\eps}} e^{-mr}\Bigl(1 +O\Bigl(\frac{m^2}{k} +
\frac{m^3}{k^2}\Bigr)\Bigr) \; \;\; + \; \; \;
o\Bigl(\frac{1}{k^2}\Bigr).
$$
To estimate the first error term we approximate the summation by an
integral, so that term becomes
 \bean
O(1) \cdot
\sum_{m=0}^{k^{1/2+\eps}}e^{-mr}  \frac{m^2}{k} & =&
O(k^{-1})\int_0^\infty
e^{-rx}x^2  \, dx\\ \\
&=&O(r^{-2}k^{-1})\\ \\
&=&O(t^{-2}k^{-1})\\ \\
&=&O\Bigl(\frac{1}{\log k}\Bigr)
 \eean
using~\eqn{tasy}. The other error term is similarly $O(1/\log k) $.
The main term in the summation is a truncated geometric series with
the truncated terms negligible, so we have
 \bean
 \pi_k &=& 1-\frac{1}{c_k}\left(o(1)+ \sum_{m=0}^{\infty}e^{-mr}\right)\\ \\
      & =& 1-\frac{1}{c_k(1-e^{-r})}+O\Bigl(\frac{1}{k\log k}\Bigr)\\ \\
      & =& 1-\frac{1}{c_k(1-(t+1)^{-1})}+O\Bigl(\frac{1}{k\log k}\Bigr)\\ \\
      & =& 1-\frac{t+1}{c_kt}+O\Bigl(\frac{1}{k\log k}\Bigr).
 \eean
Using this with~\eqn{clam1} shows that
 \bel{cktolamk}
c_k=\la_k +t^{-1} +1 +O\Bigl(\frac{1}{\log k}\Bigr).
 \ee
So we may continue with
 \bean
 \pi_k &=& 1-\frac{t+1}{\la_kt}+O\Bigl(\frac{1}{k\log k}\Bigr) \\
       &=& 1-\frac{1}{ kt}+O\Bigl(\frac{1}{k\log k}\Bigr).
 \eean
Hence
$$
\log \pi_k = -\frac{1}{ kt}+O\Bigl(\frac{1}{k\log k}\Bigr).
$$
Next, substitute this in the left side of~\eqn{logpik}, and $t=
(1+x)(q_k/k)^{1/2}$ into the right side. We know that $x=o(1)$
from~\eqn{tasy}, and we may expand $\log(1+t)$ as
$t-\frac{1}{2}t^2+\frac{1}{3}t^3+O(t^4)$. The upshot is that
$$
t = \Bigl(\frac{q_k}{k}\Bigr)^{1/2} +\frac{q_k +
3}{3k}+O\Bigl(\frac{1}{k\log k}\Bigr).
$$
This determines $t$, and since $\la_k\sim k$, we have
from~\eqn{tdef} that  $ \la_k = k+kt-2 +O(1/\log k)$. Now
using~\eqn{cktolamk} and the formula for $t$ immediately above
(which in particular gives $1/t =  (q_k/k)^{1/2} -1/3$), we obtain
Lemma~\ref{ckasy}. \qed

\bs


\begin{thebibliography}{99}

\bibitem{BCFF}
B. Bollob\'{a}s, C. Cooper, T. Fenner, and A. Frieze, Edge disjoint
Hamilton cycles in sparse random graphs of minimum degree at least
$k$. {\it J. Graph Theory} {\bf 34} (2000), no. 1, 42--59.

\bibitem{BKV} B. Bollob{\'a}s, J. H. Kim, and J. Verstra\"{e}te,
Regular subgraphs of random graphs, {\em Random Structures \&
Algorithms}, 29 (2006), 1-13.

\bibitem{bkw} I. Benjamini, G. Kozma, and N. Wormald, The mixing time of
the giant component of a random graph, {\it Preprint}.

\bibitem{CW} J. Cain and N. Wormald,  Encores on cores,
{\em  Electronic Journal of Combinatorics} {\bf 13} (2006), RP 81.

\bibitem{FKP} P. Flajolet, D. Knuth, and B. Pittel,
The first cycles in an evolving graph. Graph theory and
combinatorics (Cambridge, 1988). {\it Discrete Math.} {\bf 75}
(1989), no. 1-3, 167--215.

\bibitem{J} S. Janson, Cycles and unicyclic components in random graphs.
{\it Combin. Probab. Comput.} {\bf 12} (2003), no. 1, 27--52.

\bibitem{JL} S. Janson and M. Luczak, A simple solution to the $k$-core problem.
Tech. Report 2005:31, Uppsala.

\bibitem{JLR}
S.~Janson, T.~{\L}uczak, and A.~Ruci{\' n}ski, \emph{Random Graphs},
Wiley, New York, 2000.

\bibitem{K} J. H. Kim, Poisson cloning model for random graphs. {\it
Preprint}.

\bibitem{LP} L.~Lov\'{a}sz and M.D.~Plummer, Matching theory. North-Holland Mathematics
Studies, 121. {\it Annals of Discrete Mathematics}, 29. Budapest,
1986, xxvii+544 pp.

\bibitem{L} T. \L uczak, Size and connectivity of the $k$-core of a random graph. {\it Discrete
Math.} {\bf 91} (1991), no. 1, 61--68.


\bibitem{PSW} B. Pittel, J. Spencer, and N.  Wormald, Sudden emergence
of a giant $k$-core in a random graph, {\it J. Combinatorial Theory, Series B}
{\bf  67} (1996), 111--151.

\bibitem{PW} M. Pretti and M. Weigt, Sudden emergence of $q$-regular subgraphs in
random graphs, {\it Europhys. Lett.} 75, 8 (2006).


\end{thebibliography}
\end{document}